# Fourier Coefficients Of Some Cusp Forms
N. A. Carella

*Abstract:* The possible values of the *n*th Fourier coefficients λ(*n*) of some cusp forms of weight $k \geq 12$ are studied in this article. In particular, the values of the tau function are investigated in some details, and a potential proof of the relation $\tau(p) \neq 0$ for all primes $p \geq 2$ is provided.

Mathematics Subject Classifications: 11F25, 11F30. 11F11.

Keywords: Fourier Coefficients, Modular forms, Cusp forms, Tau functions, Lehmer conjecture.

## 1. Introduction
Let $\mathbb{H} = \{ z \in \mathbb{C} : \text{Im}(z) > 0 \}$ be the upper half plane, and let $\eta(z) = q^{1/24}\prod_{n \geq 1}(1-q^n)$ be the eta function, where $q = e^{i2\pi z}$. The eta function is holomorphic, and nonvanishing for all $z \in \mathbb{H}$. The Fourier series expansions of the powers $\eta(z)^{2k}$ of the eta function

$$\eta(z)^{2k} = q^{k/2}\prod_{n \geq 1}\left(1-q^n\right)^{2k} = \sum_{n \geq 1}\lambda(n)q^n, \tag{1.1}$$

can have vanishing Fourier coefficients $\lambda(n) = 0$ for some combinations of fixed levels $N \geq 1$, and weights $k \geq 1$ and variable $n \geq 1$. The smallest level $N = 1$, weights $k = 2, 4, 6, 8, 10, 14$ are known to have vanishing coefficients, [SE]. But, the situation for other levels $N \geq 1$ and weights $k \geq 1$ can be very different, [RT], [KR].

The discriminant function of level $N \geq 1$, and weight $k \geq 12$ is defined by $\Delta_k(z) = \eta(z)^{2k}$. The best known is discriminant function $\Delta(z) = \Delta_{12}(z)$ of level $N = 1$, and weight $k = 12$. The corresponding Fourier series is usually written as

$$\Delta(z) = \sum_{n \geq 1}\tau(n)q^n, \tag{1.2}$$

see [AP, Chapter 3], [KR, p. 20], and the *n*th coefficient $\tau(n)$ is known as the Ramanujan tau function.

The earliest works on the function $\tau : \mathbb{N} \to \mathbb{Z}$ was done by Ramanujan, Mordell, Lehmer, and other authors. The preliminary work on the question of whether or not some of the coefficients of the discriminant function $\Delta(z)$ vanish, that is, $\tau(n) = 0$ for some $n \geq 1$, was done in [LR]. Now a day more general questions are promulgated in the literature, see [RT, p. 3], [OK], et alii, for informative discussions.

*Conjecture* **1.1**. (Generalized Lehmer conjecture) Let $\Delta_k(z) = \sum_{n \geq 1}\tau_k(n)q^n \in S_k(\Gamma_0(N))$ be a new form of level $N \geq 1$, and weight $k \geq 12$. Then, $\tau_k(n) \neq 0$ for all $n \geq 1$.

It is known that the first occurrence of any vanishing coefficient occurs at prime arguments, [LR, Theorem 2].

Lehmer verified the relation $\tau_k(p) \neq 0$ for the small primes $p < 3316799$. Since then, the accumulated works of a few other authors have extended the numerical data to about $p \leq 2.210^{19}$, see [BN, p. 3], [EN, p. 170], and [SE]. A possible solution for some cusp forms of level $N = 1$ is offered in this article.

**Theorem 1.2.** Let $\Delta_k(z) = \sum_{n \geq 1} \tau_k(n) q^n \in S_k(\Gamma_0(1))$ be a cusp form of weight $k$ = 12, 16, 18, 20, 24, 26. Then, the coefficients $\tau_k(p) \neq 0$ for all primes $p \geq 2$.

The next Section 2, surveys some of the important related works in the literature. Section 3 records a few Lemmas used in the proof. The proof of Theorem 1.2, based on the congruence $\tau(n) \equiv \sigma_{11}(n) \bmod 691$, and the estimate $\#\{ p \leq x : \tau_k(p) = 0 \} = O(x(\log x)^{-3/2-\varepsilon})$, where $\varepsilon > 0$ is an arbitrarily small number, appears in Section 4, and the results for a few other values such as $\tau_k(p) \neq 2$, and $\tau_k(p) = -24$ are also given here in Sections 5 and 6. Section 7 investigates the density of the subset of integers $\{ \tau_k(n) : n \leq x \}$. Another proofof the special case $\tau_{12}(p) \neq 0$ for all primes $p \geq 2$, based on the congruence $\tau(n) \equiv n\sigma_3(n) \bmod 7$, and the estimate $\#\{ p \leq x : \tau_k(p) = 0 \} = O(x(\log x)^{-3/2-\varepsilon})$, is assembled in Theorem 8.1 in Section 8.

## 2. A Synopsis Of Previous Works

A few of the recent advances on the theory of the Fourier coefficients of the modular forms $\Delta_k(z)$ are surveyed below. These results demonstrate several interesting approaches to the theory of the Fourier coefficients of modular forms.

### 2.1 Vanishing And Lacunary Property

A function $f(z) = \sum_{n \geq 1} a(n) q^n$ is *lacunary* if the subset of nonvanishing coefficients $\{ a(n) \neq 0 : n \geq 1 \}$ has zero density in the set of integers $\mathbb{N}$. A somewhat general analysis of the lacunary properties of various powers, products, and quotients of the eta functions such as $\eta(z)^a \eta(z)^b$, $\eta(z)^a \eta(2z)^b$, etc, are presented in [KR], [OM], [SE], and similar sources.

**Theorem 2.1.** ([SE]) The cusp form $\eta^{2k}(z) \in S_k(\Gamma_0(N))$ is lacunary if and only if $k$ = 2, 4, 6, 8, 10, 14, 26.

This completes the "only if" part of the proof previously done by Ramanujan [RJ].

The case of modular forms $f(z) = \sum_{n \geq 1} a(n) q^n$ of weight $k = 2$ have vanishing coefficients $a(p) = 0$ for infinitely many primes $p$, called supersingular primes, was proved in [EL].

### 2.2 Density of the Domain

For any fixed prime $p \geq 2$, almost every $\tau_k(n)$ is divisible by $p$. This statement has the asymptotic formula $\#\{ p \leq x : \tau_k(n) \equiv 0 \bmod p \} \geq cx / \log^\beta x$, for some constants $\beta = \beta(p) > 0$, and $c > 0$.

***Theorem 2.2.*** ([SF]). Let $q$ be prime, and let $a \in \mathbb{Z}$. If there exists a prime $p \neq q$ that satisfies the congruence $\tau_k(p) \equiv a \bmod q$, then the subset of primes $\{ p \leq x : \tau_k(p) \equiv a \bmod q \}$ has positive density and has the asymptotic formula

$$\pi_k(x, q, a) = \#\{ p \leq x : \tau_k(p) \equiv a \bmod q \} \geq c(a,q) x / \log x, \tag{2.1}$$

where $c_k(a,q) > 0$ is a constant.

This is an amazing result for the existence of infinitely many primes in residue classes given that the set is not empty. For example,

$\#\{ p \leq x : \tau_{12}(p) \equiv 0 \bmod 2 \} = c_2 \pi(x)$, $\qquad \#\{ p \leq x : \tau_{12}(p) \equiv 0 \bmod 3 \} = c_3 \pi(x)$,
$\#\{ p \leq x : \tau_{12}(p) \equiv 0 \bmod 5 \} = c_5 \pi(x)$, $\qquad \#\{ p \leq x : \tau_{12}(p) \equiv 0 \bmod 7 \} = c_7 \pi(x)$,

where $c_2 = 1, c_3, c_5, c_7 > 0$ are constants.

The asymptotic formula for the density of the subset of primes (2.1) has been refined in [KM, Corollary 3.1]. Similar density results for the subsets of primes $\{ p \leq x : \tau_k(p) = m \}$ for some integers $m \in \mathbb{Z}$ are also known.

***Theorem 2.3.*** ([SF]). Given arbitrary $\varepsilon > 0$, and a large number $x \geq x_0$. The subset of primes $\{ p \leq x : \tau_k(p) = 0 \}$ has zero density in the set of prime numbers $\mathbb{P} = \{ 2, 3, 5, \ldots \}$, and satisfies the asymptotic formula

$$\#\{ p \leq x : \tau_k(p) = 0 \} = O(x/(\log x)^{3/2-\varepsilon}). \tag{2.2}$$

A sketch of the proof is given in [SF, p. 9]. A sharper estimate $\{ p \leq x : \tau_k(p) = 0 \} = O(x^{2/3})$ modulo the GRH was also proved by Serre.

***Theorem 2.4.*** Given arbitrary $\varepsilon > 0$, and a large number $x \geq x_0$, the followings hold.
(i) Let $a \geq 2$ be prime number. Then, the subset of primes $\{ p \leq x : \tau_k(p) = a \} = \emptyset$ is the empty set.

(ii) Let $a < 3316799$ be prime number. Then, the subset of primes $\{ p \leq x : \tau_k(p^2) = a \} = \emptyset$ is the empty set.

***Proof***: (i) The function $\tau_k(n)$ is odd if and only if $n \geq 1$ is a square integer. (ii) There are no integer solutions $n \in \mathbb{N}$ of the equation $\tau_k(n) = a$, for any prime $a < 3316799$, [LN]. ∎

For any fixed $a \in \mathbb{Z}$, the number of empty subsets $\{ p \leq x : \tau_k(p) = a \} = \emptyset$ is expected to be finite. And it is expected that $\#\{ p \leq x : \tau_k(p) = a \} = 1$ for almost every $a \neq 0$. A new method in Sections 5 and 6 will be used here to verify that $\#\{ p \leq x : \tau_k(p) = 2 \} = O(1)$, and $\#\{ p \leq x : \tau_k(p) = -24 \} = O(1)$.

Assuming some conjectures, in [RT] it is shown that the density of the nonzero coefficients, the support of $\tau_k(n)$, has a lower bound

$$D_k = \lim_{x \to \infty} \frac{\#\{n \leq x : \tau_k(n) \neq 0\}}{x} \geq .99999, \tag{2.3}$$

for $k = 12, 16, 18, 20, 24, 26$.

### 2.3 Density of the Range

The subset of integers $B_k = \{\tau_k(p) : \text{prime } p \geq 2\}$ is a multiplicative basis of the image $R_k = \{\tau_k(n) : \text{integer } n \geq 1\}$ of the Fourier coefficient $\tau_k(n)$. The basis has zero density in the set of integer $\mathbb{Z}$, but it is unknown if $R_k(x) \sim x$.

Unconditionally, in [GV] it is shown that the range $\tau_{12}(n)$ has a lower bound

$$R_{12}(x) = \#\{\tau_{12}(n) : n \leq x\} \geq c x^{1/2} e^{-4\log x/\log\log x}, \tag{2.4}$$

where $c > 0$ is a constant. It is believed that $R_k(x) \geq cx$ for $k = 12, 16, 18, 20, 24, 26$. The ideal case would be a density of the form $R_k(x) \sim x$ for each weight $k \geq 12$.

In [MT] it is shown that the odd values of the coefficients has a lower bound

$$|\tau_{12}(n)| \geq (\log n)^c, \tag{2.5}$$

where $c > 0$ is a constant.

### 2.4 Randomness

The values $\tau_k(p)/2p^{(k-1)/2} = \cos\theta_{k,p}$ are independent random variables over the interval $[-1, 1]$ with respect to the Sato-Tate measure $d\mu = 2\pi^{-1}\sin(\theta)^2 d\theta$. The CM case has a different distribution.

***Theorem* 2.5.** ([Deligne-Weil]). Given a prime $p \geq 2$, the $p$th coefficient of a cusp form of weight $k \geq 2$ satisfies the inequality

$$|\lambda(p)| \leq 2p^{(k-1)/2}. \tag{2.6}$$

This a generalization of the Hasse (or Deuring) inequality $|\lambda(p)| \leq 2p^{1/2}$ for cusp forms of weight $k = 2$.

## 3. Arithmetic Properties

Let $\mathbb{H} = \{z \in \mathbb{C} : \text{Im}(z) > 0\}$ be the upper half plane, and let $\text{SL}_2(\mathbb{Z}) \cong \Gamma(1)$ be the special linear group of $2 \times 2$ nonsingular integers matrices. For an integer $N \geq 1$ define the modular subgroup of level $N$ by

$$\Gamma_0(N) = \left\{ \begin{bmatrix} a & b \\ c & d \end{bmatrix} \equiv \begin{bmatrix} a & b \\ 0 & d \end{bmatrix} \bmod N \right\}. \tag{3.1}$$

The congruence subgroup $\Gamma_0(N)$ has index $[\Gamma(1):\Gamma_0(N)] = N\prod_{p|N}(1+1/p)$ in the full modular group $\Gamma(1)$.

A function $f: \mathbb{H} \to \mathbb{C}$ is a modular form of level $N \geq 1$, and weight $k \geq 1$ if

$$f\left(\frac{az+b}{cz+d}\right) = (cz+d)^k f(z) \tag{3.2}$$

for any fractional transformation $\gamma(z) = (az+b)/(cz+d)$ from the congruence subgroup $\Gamma_0(N)$. The definition immediately implies that it has a Fourier series expansion $f(z) = \sum_{n \geq 1} a(n)q^n$, where $q = e^{i2\pi z}$. Moreover, it is a cusp form if it satisfies the initial coefficients condition $a(0) = 0,\ a(1) = 1$. The eta product

$$f(z) = \prod_{d|N} \eta(dz)^{r_d}, \tag{3.3}$$

where $r_d \in \mathbb{Z}$, and $k = (1/2)\sum_{d|N} r_d$ is the weight, is a rich source of modular forms. A complete description of the eta products, and the analysis of many modular forms of various weights and levels appear in [KR], [OM, p 18], and other references in the literature.

The space of all modular forms level $N = 1$, and any weight $k \geq 4$ is denoted by $\mathcal{M} = \mathbb{C}[G_4(z), G_6(z)]$, and the subspace of all modular forms of level $N = 1$, and weight $k \geq 4$ is denoted by $\mathcal{M}_k = \mathbb{C}[G_4(z), G_6(z)]$, where $G_k(z) = \sum_{(m,n) \neq (0,0)} (mz+n)^{-k}$, $z \in \mathbb{H}$, is the Eisenstein series of weight $k \geq 4$, see [IH, p. 20].

***Theorem* 3.1.** The subspace $\mathcal{M}_k$ of modular forms of level $N = 1$, and weight $k \geq 4$ has dimension

$$\dim M_k = \begin{cases} [k/12] & \text{if } k \equiv 2 \bmod 12, \\ [k/12]+1 & \text{if } k \not\equiv 2 \bmod 12. \end{cases} \tag{3.4}$$

**3.1 Coefficients Congruences**

The trace $Tr_{2k}(\Gamma_0(N), n) = \lambda(n)$ of the linear Hecke operator acting on the set of function is defined by

$$T_n f(z) = \sum_{m \geq 1,\ d|\gcd(m,n),\ \gcd(d,N)=1} d^{2k-1}\lambda(mn/d^2)\, q^m = \lambda(n). \tag{3.5}$$

The *n*th Fourier coefficients $\lambda(n)$ of modular forms $f(z) = \sum_{n \geq 1} \lambda(n)q^n \in S_k(\Gamma_0(N))$ are multiplicative, and satisfy an amazing array of congruences and recurring formulas. A few of these relations, which are used in the proofs, are stated here.

Fourier Coefficients Of Some Cusp Forms

In the following the notation $\sigma_s(n) = \sum_{d \mid n} d^s$ is the usual sum of divisors function. And a prime $q \neq p$ is exceptional if the image of the representation ρ is not all of $GL_2(\mathbb{Z}_p)$, [LA, p. 196], [BN], et alii for further explanation.

***Theorem 3.2.*** Let $\tau_k(n)$ be the $n$th coefficient of the cusp form $\Delta_k(z) = \sum_{n \geq 1} \tau_k(n) q^n \in S_k(\Gamma_0(1))$. Then

(i) $\tau_{12}(n) \equiv \sigma_{11}(n) \bmod 691$.   The exceptional primes are 2, 3, 5, 7, 23, and 691.

(ii) $\tau_{16}(n) \equiv \sigma_{15}(n) \bmod 3617$.   The exceptional primes are 2, 3, 5, 7, 11, 31, 59?, and 3617.

(iii) $\tau_{18}(n) \equiv \sigma_{17}(n) \bmod 43867$.   The exceptional primes are 2, 3, 5, 7, 11, 13, and 43867.

(iv) $\tau_{20}(n) \equiv \sigma_{19}(n) \bmod 617$.   The exceptional primes are 2, 3, 5, 7, 11, 13, 283, and 617.

(v) $\tau_{24}(n) \equiv \sigma_{23}(n) \bmod 593$.   The exceptional primes are 2, 3, 5, 7, 13, 131, and 593.

(vi) $\tau_{26}(n) \equiv \sigma_{25}(n) \bmod 657931$.   The exceptional primes are 2, 3, 5, 7, 11, 17, and 657931.

The development of these congruences, and other formulae are given in [SW], [RT], and [BT]. The proof of congruence (i), in various levels of details, appears in [AP], [LM], [BT], [SF, p. 4], [ZR], and many other sources. The nonordinary primes $p$ such that $\tau_k(p) \equiv 0 \bmod p$ re studied in [GF], and [LZ].

***Lemma 3.3.*** (Mordell) For every prime $p \geq 2$, and $n \geq 1$, the $n$th coefficient $\lambda(n)$ of a modular form $f(z) = \sum_{n \geq 1} \lambda(n) q^n \in S_k(\Gamma_0(N))$ of level $N \geq 1$, and weight $k \geq 2$ satisfy the multiplication formulas

(i) $\lambda(mn) = \lambda(m)\lambda(n)$,   $\gcd(m, n) = 1$.   (3.6)
(ii) $\lambda(p^{n+1}) = \lambda(p)\lambda(p^n) - p^{k-1}\lambda(p^{n-1})$,   $p \geq 2$ prime and $n \geq 1$.
(iii) $\lambda(mp) = \lambda(m)\lambda(p)$,   for prime $p \mid N$.

The proof assumes that the Dirichlet series of the modular function has a product expansion of the shape

$$L(s, f) = \sum_{n \geq 1} \frac{\lambda(n)}{n^s} = \prod_{\gcd(N,p) \neq 1} \left(1 - \lambda(p) p^{-s}\right)^{-1} \prod_{\gcd(N,p)=1} \left(1 - \lambda(p) p^{-s} + p^{k-1-2s}\right)^{-1},$$   (3.7)

refer to [AP, p. 92], [LA], etc.

***Lemma 3.4.*** For every prime $p \geq 2$, and $n \geq 1$, the $p$th coefficient $\lambda(p)$ of a modular form $f(z)$ satisfy the congruence $\lambda(p^n) \equiv \lambda(p)^n \bmod p$.

***Proof***: Use induction on the multiplicative formula.   ∎

***Theorem* 3.5.** (Hecke) If $f(z) = \sum_{n \geq 1} \lambda(n) q^n$ is a modular forms of level $N \geq 1$, and weight $k \geq 2$, then the L-function $L(s, f)$ has an analytic continuation to the entire complex plane $\mathbb{C}$.

The L-function $L(s, f)$ has another description as the Mellin transform of the modular function.

***Theorem* 3.6.** If $f(z) = \sum_{n \geq 1} \lambda(n) q^n$ is a modular forms of level $N \geq 1$, and weight $k \geq 2$, then

$$L(s, f) = \frac{(2\pi)^s}{\Gamma(s)} \int_0^\infty f(iy) y^{s-1} \, dy, \tag{3.8}$$

where $z = x + iy \in \mathbb{H}$.

## 3.2 Algorithms

The algorithmic aspects of modular functions are cover in extensive details in [EN], and [ST]. The analysis and other information on the algorithm for computing large values $\tau(p)$, $p \leq 10^{13}$, are given in [LZ]. The time complexity of this algorithm is developed in [EN], see also [ST, p. 29].

***Lemma* 3.7.** ([EN]) For every prime $p \geq 2$, the $p$th coefficient $\lambda(p)$ of a modular form $f(z) = \sum_{n \geq 1} \lambda(n) q^n$ of weight $k \geq 2$ can be computed in polynomial time $O((\log p)^c)$, $c > 0$ constant.

## 3.3 A Result On Diophantine Equations

A few authors have worked on the quantitative version of the qualitative results of Thue and Siegel for the finiteness of integral points on algebraic curves. A recent quantized version states the followings.

***Theorem* 3.8.** ([SN]) Let $K$ be a numbers field, let $R_S$ be the ring of $S$-integers of $K$ for some finite set of places $S$, and let $E : y^2 = x^3 + ax + b$ be an elliptic curve in quasi-minimal Weierstrass equation. Then

(i) If the $j$-invariant is integral, there is a constant $\kappa > 1$, depending on $K$ such that

$$\#\{ P \in E(K) : x(P) \in R_S \} \leq \kappa^{1 + \#S + \mathrm{rank}(E(K))}, \tag{3.9}$$

where $\mathrm{rank}(E(K)) \geq 0$ is the rank of $E$.

(ii) If the $j$-invariant is nonintegral for at most $\delta \geq 1$ primes of $K$, there is a constant $\kappa > 1$, depending on $K$ such that

$$\#\{ P \in E(K) : x(P) \in R_S \} \leq \kappa^{1 + \#S + (1+\delta)\mathrm{rank}(E(K))}, \tag{3.10}$$

where $\mathrm{rank}(E(K)) \geq 0$ is the rank of $E$.

A quasi-minimal Weierstrass equation $E : y^2 = x^3 + ax + b$ has minimal discriminant $\left| N_{K/\mathbb{Q}}(4a^3 + 27b^2) \right|$

subject to integral $a, b \in K$.

### 3.4 Large Integral Solutions

Various results for the maximal size of the largest integral solutions of certain Diophantine equations have been achieved by a few workers in the field of Diophantine analysis. The theory is covered in fine details in [SV, Chapter 6].

The height of a polynomial $f(x, y) = \sum_{0 \leq i, j \leq d} a_{i,j} x^i y^j$ of degree $\deg f(x, y) = d$ is defined by $H = \max \{ |a_{i,j}| : 0 \leq i, j \leq d \}$. And the discriminant is defined by the product $\prod_{i \neq j} (\alpha_i - \alpha_j)^2$ over the roots of the polynomial $f(\alpha, 1) = 0$.

***Theorem* 3.9.** ([SV, p. 115]) Let $f(x, y) = 0$ be a cubic polynomial of nonzero discriminant, and height $H$. Then

$$\max \{ |x|, |y| : f(x, y) = 0 \} < e^{cH^{6+\varepsilon}}, \tag{3.11}$$

where $c > 0$ is a constant, and $\varepsilon > 0$ is arbitrarily small.

### 3.5 An Estimate Of The Rank

The height $H = \max \{ |a_i| : 0 \leq i \leq 6 \}$ of an elliptic curve $E : y^2 + a_1 xy + a_3 y = x^3 + a_2 x^2 + a_4 x + a_6$ over the integers $\mathbb{Z}$ restricts the size of the rank $\operatorname{rank}(E) \geq 0$ of the group of rational points $E(\mathbb{Q}) \cong \mathbb{Z}^{\operatorname{rank}(E)} \times E(\mathbb{Q})_{\operatorname{tor}}$. Elliptic curves with groups of rational points of large ranks have large coefficients, but not conversely.

The Birch-Swinnerton-Dyer conjecture states that

$$\pi_E(x) = \prod_{p \leq x} \frac{N_p}{p} \sim c(\log x)^{\operatorname{rank}(E)}, \tag{3.12}$$

where $p \nmid \Delta(E)$, which is the discriminant of the elliptic curve, $N_p = \#E(\mathbf{F}_p)$, and $c > 0$ is a constant, see [RS]. The cardinality of the group of points $E(\mathbf{F}_p) = \{ (x, y) : f(x, y) \equiv 0 \bmod p \}$ is in the range $p + 1 - 2\sqrt{p} \leq N_p \leq p + 1 + 2\sqrt{p}$. In fact, the number $N_p = p + 1 + 2\sqrt{p} \cos \theta_p$ is an independent random variable with respect to the Sato-Tate measure $d\mu = 2\pi^{-1} \sin(\theta)^2 d\theta$. Hence,

$$\prod_{p \leq X} \left( 1 + \frac{1}{p} + \frac{2 \cos \theta_p}{\sqrt{p}} \right) \sim c(\log X)^{\operatorname{rank}(E)}. \tag{3.13}$$

This has the equivalent asymptotic $\log(c \log X)^{-1} \sum_{p \leq X} \log\left(1 + p^{-1} + 2p^{-1/2} \cos \theta_p\right) \sim \operatorname{rank}(E)$, certain finite sums related to this are studied in [FD]. On average, an elliptic curve is expected to have a rank of 1/2, (this a conjecture of Goldfeld).

***Lemma* 3.10.** Let $E: f(x, y) = 0$ be an elliptic curve of height $H^c \leq X$, $c > 2$ constant, then

$$\text{rank}(E) \leq c_0 \log X / \log \log X, \tag{3.14}$$

where $c_0 > 0$ is a constant.

A sharper upper bound for the case $E: y^2 = x^3 + a_2 x^2 + a_4 x$ is given in [AU].

### 3.6 The Fourier Series of the Discriminant Functions

Each $\Delta_k(z) = \sum_{n \geq 1} \tau_k(n) q^n$ generates a one dimensional space of cusp forms $S_k(\Gamma_0(1))$, $k = 12, 16, 18, 20, 24, 26$. The series expansions of the discriminant functions $\Delta_k(z) = \sum_{n \geq 1} \tau_k(n) q^n \in S_k(\Gamma_0(1))$, [AP, p. 133], have the following formulas specified in terms of the first cusp form, and the Eisenstein series:

(i) $\Delta_{12}(z) = \sum_{n \geq 1} \tau_{12}(n) q^n = q - 24q^2 + 252q^3 - 1472q^4 + 4830q^5 - 6048q^6 - 16744q^7 + 84480q^8 + \cdots$.

(3.14)

(ii) $\Delta_k(z) = \dfrac{\Delta_{12}(z)}{(2\pi)^2} E_{2k-12}(z) = \sum_{n \geq 1} \tau_{12}(n) q^n \left( 1 - \dfrac{2(2k-12)}{B_{2k-12}} \sum_{m \geq 1} \sigma_{2k-13}(m) q^m \right)$ (3.15)

for $k = 16, 18, 20, 24, 26$. The $n$th coefficient has the convenient formula

$$\tau_k(z) = -\dfrac{4k - 24}{B_{2k-12}} \sum_{0 \leq m \leq n} \tau_{12}(m) \sigma_{2k-13}(n - m). \tag{3.16}$$

The Eisentein series is given by

$$E_{2k}(z) = 1 - \dfrac{2 \cdot 2k}{B_{2k}} \sum_{m \geq 1} \sigma_{2k-1}(m) q^m \tag{3.17}$$

for $k \geq 1$, [AP, p. 139], [KZ, p. 111]. The first six are these:

$$E_4(z) = 1 + 240 \sum_{n \geq 1} \sigma_3(n) q^n, \qquad E_6(z) = 1 - 504 \sum_{n \geq 1} \sigma_5(n) q^n, \tag{3.18}$$

$$E_8(z) = 1 + 480 \sum_{n \geq 1} \sigma_7(n) q^n, \qquad E_{10}(z) = 1 - 264 \sum_{n \geq 1} \sigma_9(n) q^n,$$

$$E_{12}(z) = 1 + \dfrac{65520}{691} \sum_{n \geq 1} \sigma_{11}(n) q^n, \qquad E_{14}(z) = 1 - 24 \sum_{n \geq 1} \sigma_{13}(n) q^n.$$

The $n$th Bernoulli number $B_n$ is defined by

$$\sum_{n \geq 0} B_n \frac{x^n}{n!} = 1 - \frac{1}{2}\frac{x}{1!} + \frac{1}{6}\frac{x^2}{2!} - \frac{1}{30}\frac{x^4}{4!} + \frac{1}{42}\frac{x^6}{6!} - \frac{1}{30}\frac{x^8}{8!} + \frac{1}{65}\frac{x^{10}}{10!} - \cdots. \tag{3.19}$$

The discriminant cusp form $\Delta_{12}(z) = \sum_{n \geq 1} \tau_{12}(n) q^n$ has the corresponding Dirichlet series

$$L(s, \Delta_k) = \sum_{n \geq 1} \frac{\tau_k(n)}{n^s} = \prod_{p \geq 2} \left(1 - \tau_k(p) p^{-s} + p^{k-1-2s}\right)^{-1}. \tag{3.20}$$

This series has an analytic continuation to an entire function, and satisfies the functional equation

$$(2\pi)^{-s} \Gamma(s) L(s, f) = (-1)^{k/2} (2\pi)^{k-s} \Gamma(k-s) L(k-s, f). \tag{3.21}$$

**3.4 Examples of Lacunary Eta Products**
The first few eta products are well known to be lacunary Fourier series.

$k = 1$. The first one is the Euler pentagonal numbers series:

$$\eta(z) = \prod_{n \geq 1} (1 - q^n) = \sum_{n \geq -\infty} (-1)^n q^{(3n^2 + n)/2}. \tag{3.22}$$

$k = 2$. The second lacunary eta product is related to the modular form $f(z) = \sum_{n \geq 1} a(n) q^n \in S_k(\Gamma_0(N))$ attached to an elliptic curve. The $p$th coefficient vanish $a(p) = 0$ at each supersingular prime $p$. The corresponding Dirichlet series

$$L(s, f) = \sum_{n \geq 1} \frac{a(n)}{n^s} = \prod_{\gcd(N, p) \neq 1} \left(1 - \lambda(p) p^{-s}\right)^{-1} \prod_{\gcd(N, p) = 1} \left(1 - \lambda(p) p^{-s} + p^{k-1-2s}\right)^{-1}, \tag{3.23}$$

has an analytic continuation to an entire function, and satisfies the functional equation

$$\left(N^{1/2} / 2\pi\right)^s \Gamma(s) L(s, f) = \omega \left(N^{1/2} / 2\pi\right)^{2-s} \Gamma(2-s) L(2-s, f), \tag{3.24}$$

where $N$ is the conductor, and $\omega = \pm 1$ is the root number, [IH, p. 136].

$k = 3$. The third one is the Jacobi triple product:

$$\eta(z)^3 = \prod_{n \geq 1} (1 - q^n)^3 = \sum_{n \geq 0} (-1)^n (2n+1) q^{(n^2 + n)/2}. \tag{3.25}$$

## 4. The Values Of the Fourier Coefficients

The congruences and recursive formulas satisfied by the *n*th coefficient $\lambda(n)$ of a modular form $f(z) = \sum_{n \geq 1} \lambda(n) q^n \in S_k(\Gamma_0(N))$ restrict the image $\lambda(\mathbb{N}) = \{ \lambda(n) \in \mathbb{Z} : n \geq 1 \}$ of the integer-valued function $\lambda : \mathbb{N} \to \mathbb{Z}$ to a subset of the integers $\lambda(\mathbb{N}) \subset \mathbb{Z}$, depending on the level $N \geq 1$ and the weight $k \geq 1$. Likely, the function $\lambda$ is injective, and the image $\lambda(\mathbb{N})$ is a proper subset of the integers $\mathbb{Z}$.

In the case of the tau function $\tau : \mathbb{N} \to \mathbb{Z}$, an odd value $\tau(n) = 2m+1$ occurs if and only if $n = d^2$, where $d \geq 1$ is an odd integer. This is derived from a congruence $\tau(n) \equiv \sigma_{11}(n) \bmod 2^v$, or by other means, see Theorem 2, [BT], [MT], et cetera. And a prime value $\tau(n) = r$ occurs if and only if $n = p^{q-1}$, where $p$ and $q \geq 3$ are primes. This is derived from the divisibility property of the quadratic recurrence relation satisfied by $\tau(n)$, see [LS].

Some criteria in terms of Diophantine equations for the existence certain integer values $m \in \mathbb{Z}$ of the coefficients $\tau_k(n) = m$ are sketched here.

Given an arbitrary fixed integer $m \in \mathbb{Z}$, the system of Diophantine equations

$$\begin{aligned} \tau_k(a_1 n_1) &\equiv \sigma_{k-1}(a_1 n_1) \bmod q_1, \\ \tau_k(a_2 n_2) &\equiv \sigma_{k-1}(a_2 n_2) \bmod q_2, \\ \tau_k(a_3 n_3) &\equiv \sigma_{k-1}(a_3 n_3) \bmod q_3, \end{aligned} \quad (4.1)$$

where $a_i, n_i \in \mathbb{N}$ are suitable integers, and the $q_i$ are prime powers, should provide sufficient information to determine if the *n*th Fourier coefficient satisfies the equation $\tau_k(n) = m$ for some $n \in \mathbb{N}$. Several challenging problems are encountered while answering this question.

(1) What is the prime factorization of *n*? This could be a difficult problem due to the lack of the completely multiplicative property of $\tau_k(n)$. For example, it is known that $\tau_k(n) = 7$ can occur if and only if $n = p^{q-1}$, where $p$ and $q \geq 3$ are primes. But what can be said about the form of the integer *n* for which $\tau_k(n) = 7^4$?

(2) The calculations of the integral points of the associated system of Diophantine equations could be a difficult problem.

This key idea is employed here to probe the existence of a few of the unknown values for prime argument $p \geq 2$, for example,

(i) $\tau_k(n) = 0$, (ii) $\tau_k(n) = \pm 1$, (iii) $\tau_k(n) = \pm 2$, etc.

### 4.1 The Value $\tau_k(n) = 0$

To study the integer solutions $n \in \mathbb{N}$ of the equation $\tau_k(n) = 0$, it is sufficient to consider prime argument $n = p \geq 2$, this arises from the multiplicative property $\tau(pq) = \tau(p)\tau(q)$, $\gcd(p,q) = 1$ of the coefficients, see [LR, Theorem 2]. The verification of whether or not $\tau(p) = 0$ for small primes $p \geq 2$ can be done by hands or machine calculations, and congruences. The verification of $\tau(p) \neq 0$ for $p \leq 2.2 \cdot 10^{19}$ is undertaken in [BN, p.

3]. The numerical values of $\tau(p) \neq 0$ for $p \leq 10^9$, are listed in [LS]. And for a list of some of the values $\tau_k(p)$ for $p < 10^6$, and $k = 12, 16, 18, 20, 22, 26$, see [GF].

**Theorem 1.2.** Let $\Delta_k(z) = \sum_{n \geq 1} \tau_k(n) q^n \in S_k(\Gamma_0(1))$ be a cusp form of weight $k = 12, 16, 18, 20, 24, 26$. Then, the coefficients $\tau_k(p) \neq 0$ for all primes $p \geq 2$.

**Proof**: Without loss in generality, let $\tau(p) = \tau_{12}(p)$, and suppose that $\tau(p) = 0$ for some prime $p \neq 691$. Now consider the congruence

$$\tau(p) \equiv \sigma_{11}(p) \bmod 691, \tag{4.2}$$

see Theorem 3.2. Expanding the divisors function $\sigma_s(n) = \sum_{d \mid n} d^s$ produces the Diophantine equation

$$p^{11} + 1 + 691v = \tau(p), \tag{4.3}$$

where $v \in \mathbb{Z}$ is a rational integer, and the corresponding congruence equation

$$p^{11} + 1 + 691v \equiv \tau(p) \bmod 691. \tag{4.4}$$

The hypothesis $\tau(p) = 0$ implies that $\tau(p) \equiv 0 \bmod 691$, and by Theorem 2.2, the subset of primes

$$\{ p \leq x : \tau(p) \equiv 0 \bmod 691 \} \tag{4.5}$$

has positive density. Therefore, equation (4.4) has infinitely many integer solutions. In fact, an explicit sequence of solutions is given by

$$v = -(p^{11} + 1)/691, \tag{4.6}$$

where the prime $p$ ranges over some arithmetic progressions, for example, $p = 2 \cdot 691m - 1$, $m \geq 1$. This sequence of solutions is the only possibility for both (4.3) and (4.4); other solutions of (4.4) yield $\tau(p) \neq 0$.

This in turns implies that the congruence $\tau(p) \equiv 0 \bmod 691$ holds over the integers, that is, $\tau(p) = 0$ for each integer solution of the form (4.6). In particular, the subset of primes

$$\{ p = 2 \cdot 691m - 1 \leq x : \tau(p) = 0 \}, \tag{4.7}$$

has positive density in the set of primes. But, this contradicts the zero density result for the subset of primes $\{ p \leq x : \tau(p) = 0 \}$, see Theorem 2.3. Specifically,

$$c_0 \frac{x}{\log x} \leq \#\{ p = 2 \cdot 691m - 1 \leq x : \tau(p) = 0 \}$$
$$\leq \#\{ p \leq x : \tau(p) = 0 \} \qquad (4.8)$$
$$\leq c_1 \frac{x}{(\log x)^{3/2-\varepsilon}},$$

where $c_0, c_1 > 0$ are constants, and $\varepsilon > 0$ is arbitrarily small. Ergo, $\tau(p) \neq 0$ for all primes $p \geq 2$.

The analysis for the Fourier coefficients $\tau_{12}(n)$ of the cusp form $\Delta_{12}(z) = \sum_{n \geq 1} \tau_{12}(n) q^n$ is applicable to the other cusp forms $\Delta_k(z) = \sum_{n \geq 1} \tau_k(n) q^n$ of weight $k = 12, 16, 18, 20, 24, 26$, mutatis mutandis. There are minor variations on the coefficients of the equations depending on the weight $k$, and the corresponding congruence in Theorem 2.3. ∎

As the prime $p \geq 2$ varies over the set of primes, the $p$th coefficient $\tau(p)$ is in the range $-2p^{11/2} \leq \tau(p) \leq 2p^{11/2}$, this forces every solution of the equation (4.3) to be a large negative integer in the range

$$\frac{-p^{11} - 1 - 2p^{11/2}}{691} < v < \frac{-p^{11} - 1 + 2p^{11/2}}{691}, \qquad v \neq \frac{-p^{11} - 1}{691}. \qquad (4.9)$$

This follows from the inequality $0 < |\tau(p)| < 2p^{(k-1)/2}$ for $k = 12$, see Theorem 2.5.

In the case of cusp forms $f(z) = \sum_{n \geq 1} a(n) q^n$ of level $N \geq 1$, and weight $k = 2$, the $p$th coefficient $a_p$ is in the range $-2p^{1/2} \leq a_p \leq 2p^{1/2}$, and the value $a_p = 0$ occurs infinitely often, so $0 \leq |a(p)| \leq 2p^{1/2}$. There seems to be no congruence of the form $a(n) \equiv \sigma_{k-1}(n) \bmod q$, $q$ a fixed prime, nevertheless, it has a simple proof. But the proofs for the frequencies of other values $a_p = a \neq 0$ are not simple to prove, in fact, this is an open problem known as Lang-Trotter conjecture.

An interesting question concerns the signs changes of the $p$th coefficients $\tau_k(p)$ of a cusp forms of large weight $k \geq 12$ for primes in arithmetic progressions. What are the asymptotic formulae:

$$T_+(x) = \#\{ \tau(p) > 0 : p \leq x, p \equiv a \bmod q \} \stackrel{?}{=} \pi(x)/2\varphi(q), \qquad (4.9)$$

and

$$T_-(x) = \#\{ \tau(p) < 0 : p \leq x, p \equiv a \bmod q \} \stackrel{?}{=} \pi(x)/2\varphi(q),$$

for $0 < a < q$, such that $\gcd(a, q) = 1$?

## 5 The Value $\tau_k(n) = 2$

There are no integer solutions $n \in \mathbb{N}$ of the equation $\tau_k(n) = p$, for any prime $p < 3316799$, [LR]. The new method will be used here to verifies this for $p = 2$ for all sufficiently large $n \geq 1$. Any solution of $\tau_k(n) = 2$ requires a prime $n \in \mathbb{N}$ argument. And since the Fourier coefficients are not completely multiplicative, the existence or nonexistence of such solution seems to imply the same fate for $\tau_k(p) = 2^m$, $m \geq 1$.

**Theorem 5.1.** Let $\Delta_k(z) = \sum_{n \geq 1} \tau_k(n) q^n \in S_k(\Gamma_0(1))$ be a cusp form of weight $k = 12, 16, 18, 20, 24, 26$. Then, the coefficients $\tau_k(n) \neq 2$ for all but finitely many integers $n \geq 1$.

**Proof**: Without loss in generality, let $\tau(p) = \tau_{12}(p)$, and suppose that $\tau(p) = 2$ for some prime $p > 2$. Consider the system of congruences

$$\begin{aligned} \tau(p) &\equiv \sigma_{11}(p) \bmod 691, \\ \tau(p^2) &\equiv \sigma_{11}(p^2) \bmod 691, \\ \tau(p^3) &\equiv \sigma_{11}(p^3) \bmod 691, \end{aligned} \tag{5.1}$$

see Theorem 3.2. Expanding $\sigma_s(n) = \sum_{d \mid n} d^s$, and transforming (5.1) into a system of Diophantine equations produces

$$\begin{aligned} \tau(p) &= p^{11} + 1 + 691v, \\ \tau(p^2) &= p^{22} + p^{11} + 1 + 691x, \\ \tau(p^3) &= p^{33} + p^{22} + p^{11} + 1 + 691y, \end{aligned} \tag{5.2}$$

where $v, x, y \in \mathbb{Z}$ are rational integers. By hypothesis $\tau(p) = 2$, and from Lemma 3.3, it quickly follows that

$$\tau(p^2) = 2^2 - p^{11} \quad \text{and} \quad \tau(p^3) = 2(2^2 - p^{11}) - p^{11} \cdot 2. \tag{5.3}$$

Substituting these values, and $u = p^{11}$ into the penultimate system of equations returns

$$\begin{aligned} u + 1 + 691v &= 2, \\ u^2 + u + 1 + 691x &= 2^2 - u, \\ u^3 + u^2 + u + 1 + 691y &= 2(2^2 - u) - 2 \cdot u. \end{aligned} \tag{5.4}$$

Eliminating $u$ and $v$, by algebraic means, Grobner basis, or resultant calculation, etc, yields

$$477481 x^3 - 20 \cdot 691 x^2 + 100 x + 2764 xy - 40 y + 691 y^2 = 0. \tag{5.5}$$

The change of variables $(x, y) \to (-X/691, Y/691)$ in (5.5) produces an elliptic curve in Weierstrass form

$$Y^2 - 4XY - 40Y = X^3 + 20X^2 + 4X. \tag{5.6}$$

Since the discriminant and the *j*-invariant

$$\Delta(E) = -2^{15} \cdot 3^2 \cdot 5 \cdot 7 \cdot 17 \quad \text{and} \quad j(E) = -472392/595, \tag{5.7}$$

respectively, are nonvanishing, the algebraic curve is nonsingular over the rational numbers $\mathbb{Q}$.

By Theorem 6 there are a finite number of integral solutions. Some of the integer points of (5.5) are:
(1) $(0, 0)$,
(2) $(0, 40/691)$,
(3) $(-687, 474727)$,
(4) $(-687, 326137449/691)$,
(5) $(-695, 480255)$,
(6) $(-695, 333777225/691)$, ... .

Replacing the third point $(x, y) = (-687, 474727)$ in (5.4) gives a nonprime power solution $u = 2^4 \cdot 43$, and the fifth point $(x, y) = (-695, 480255)$ in (5.4) gives a nonprime power solution $u = 2^2 \cdot 173$, ... . ∎

To settle the finite number of potential primes $p$ such that $\tau(p) = 2$, the integral solutions of (5.5) or (5.6) are required. Equations (5.5) and (5.6) are rationally equivalent, so these algebraic equations have the same group of rational points $C(\mathbb{Q}) \cong E(\mathbb{Q})$ up a finite subset of points. The rank of the elliptic curve (5.6) is rank $E(\mathbb{Q}) \geq 1$, and the presence of large rational coefficients, perhaps, suggest that it is unlikely to have other integers points.

By Theorem 3.8, the algebraic curve (5.5) or equivalently (5.6) has finitely many integer solutions $(x \leq X_0, y \leq Y_0)$, where $X_0 = X_0(k)$ and $Y_0 = Y_0(k)$ are constants depending on $k$, see Theorem 3.9.

## 6 The Value $\tau_k(n) = -24$

The prime factorization of any solution $n \in \mathbb{N}$ of $\tau_k(n) = -24$ is not clear at this point. It will be assumed it is a prime $n = p \geq 2$ argument. This assumption facilitates the calculation of $\sigma_s(n) = \sum_{d \mid n} d^s$, see below.

**Theorem 6.1.** Let $\Delta_k(z) = \sum_{n \geq 1} \tau_k(n) q^n \in S_k(\Gamma_0(1))$ be a cusp form of weight $k = 12, 16, 18, 20, 24, 26$. Then, the coefficients $\tau_k(n) = -24$ has at most a finite number of integer solutions $n \geq 1$.

**Proof**: Without loss in generality, let $\tau(p) = \tau_{12}(p)$, and suppose that $\tau(p) = -24$ for some prime $p \geq 2$. Consider the system of congruences

$$\begin{aligned} \tau(p) &\equiv \sigma_{11}(p) \bmod 691, \\ \tau(p^2) &\equiv \sigma_{11}(p^2) \bmod 691, \\ \tau(p^3) &\equiv \sigma_{11}(p^3) \bmod 691, \end{aligned} \tag{6.1}$$

see Theorem 3.2. Expanding $\sigma_s(n) = \sum_{d \mid n} d^s$, and transforming (6.1) into a system of Diophantine equations produces

$$\begin{aligned}\tau(p) &= p^{11} + 1 + 691v, \\ \tau(p^2) &= p^{22} + p^{11} + 1 + 691x, \\ \tau(p^3) &= p^{33} + p^{22} + p^{11} + 1 + 691y,\end{aligned} \qquad (6.2)$$

where $v, x, y \in \mathbb{Z}$ are rational integers.

By the hypothesis $\tau(p) = -24$, and from Lemma 3.3, it quickly follows that

$$\tau(p^2) = (-24)^2 - p^{11} \quad \text{and} \quad \tau(p^3) = (-24)((-24)^2 - p^{11}) - p^{11} \cdot (-24). \qquad (6.3)$$

Substituting these values, and $u = p^{11}$ into the penultimate system of equations returns

$$\begin{aligned}u + 1 + 691v &= -24, \\ u^2 + u + 1 + 691x &= 24^2 - u, \\ u^3 + u^2 + u + 1 + 691y &= -24(24^2 - u) + 24u.\end{aligned} \qquad (6.4)$$

Eliminating $u$ and $v$, by algebraic means, Grobner basis, or resultant calculation, etc, yields

$$477481x^3 - 1127712x^2 + 942340x + 2764xy + 25440y + 691y^2 = 0. \qquad (6.5)$$

The change of variables $(x, y) \to (-X/691, Y/691)$ in (6.5) produces an elliptic curve in Weierstrass form

$$Y^2 + 4XY + 25440Y = X^3 + 1632X^2 + 942340X. \qquad (6.6)$$

Since the discriminant and the $j$-invariant

$$\Delta(E) = -2^{17} \cdot 3^2 \cdot 5^3 \cdot 23 \cdot 53^3 \cdot 19031 \quad \text{and} \quad j(E) = -870728796277342/73311073088625, \qquad (6.7)$$

respectively, are nonvanishing, the algebraic curve is nonsingular over the rational numbers $\mathbb{Q}$.
By Theorem 3.6 there are a finite number of integral solutions. Some of the integer points of (6.5) are:

(1) $(0, 0)$,  (2) $(0, 25440/691)$,
(3) $(-675, -12437115)$,  (4) $(-675, 8610812325/691)$, ... .

Replacing the third point $(x, y) = (-675, -12437115)$ in (6.4) gives a prime power solution $u = 2^{11}$. This confirms the known value $\tau_{12}(2) = -24$. Moreover, since the curve is of genus 1, there are at most finite many more integer solutions $(u = p^{11}, v, x, y)$ of the algebraic equation (6.5) such that $\tau_{12}(p) = -24$. ∎

To settle the finite number of potential primes $p$ such that $\tau(p) = 2$, the integral solutions of (6.5) or (6.6) are required. Equations (6.5) and (6.6) are rationally equivalent, so these algebraic equations have the same group of rational points $C(\mathbb{Q}) \cong E(\mathbb{Q})$ up a finite subset of points. The rank of the elliptic curve (6.5) is rank $E(\mathbb{Q}) \geq$ 1, and the presence of large rational coefficients, perhaps, suggest that it is unlikely to have other integers points.

By Theorem 3.8, the algebraic curve (6.5) or equivalently (6.6) has finitely many integer solutions $(x \leq X_0, y \leq Y_0)$, where $X_0 = X_0(k)$ and $Y_0 = Y_0(k)$ are constants depending on $k$, see Theorem 3.9.

# 7 Density Of The Values $\tau_k(n)$

There is some evidence to suggest that the map $n \to \tau_k(n)$ is injective. The current literature claims that the density of the values of the Fourier coefficients of $\Delta_{12}(z) = \sum_{n \geq 1} \tau_{12}(n) q^n \in S_{12}(\Gamma_0(1))$ has the lower bound $\#\{\tau_{12}(n) : n \leq x\} \geq cx^{1/2} e^{-4\log x / \log\log x}$, $c > 0$ constant, see [GV]. The technique used in the previous Sections are applicable to this problem too.

**Theorem 7.1.** Let $\Delta_k(z) = \sum_{n \geq 1} \tau_k(n) q^n \in S_k(\Gamma_0(1))$ be a cusp form of weight $k$ = 12, 16, 18, 20, 24, 26. Then, density of the values of the coefficients satisfies the lower bound.

$$\#\{\tau_k(n) : n \leq x\} \geq c_k x^{1+o(1)} / \log x, \tag{7.1}$$

where $c_k > 0$ is a constant depending on the weight $k$.

***Proof***: Without loss in generality, let $\tau(p) = \tau_{12}(p)$. Let $x \geq x_0$ be a large number, and fix an integer $t \in \mathbb{Z}$ such that $|t| \leq x^{11/2+\varepsilon}$ and suppose that $\tau(p) = t$. Consider the system of congruences

$$\begin{aligned} \tau(p) &\equiv \sigma_{11}(p) \bmod 691, \\ \tau(p^2) &\equiv \sigma_{11}(p^2) \bmod 691, \\ \tau(p^3) &\equiv \sigma_{11}(p^3) \bmod 691, \end{aligned} \tag{7.2}$$

see Theorem 3.2. Expanding $\sigma_s(n) = \sum_{d \mid n} d^s$, and transforming (7.2) into a system of Diophantine equations produces

$$\begin{aligned} \tau(p) &= p^{11} + 1 + 691v, \\ \tau(p^2) &= p^{22} + p^{11} + 1 + 691x, \\ \tau(p^3) &= p^{33} + p^{22} + p^{11} + 1 + 691y, \end{aligned} \tag{7.3}$$

where $v, x, y \in \mathbb{Z}$ are rational integers.

By the hypothesis $\tau(p) = t$, and from Lemma 3.3, it quickly follows that

$$\tau(p^2) = t^2 - p^{11} \quad \text{and} \quad \tau(p^3) = t(t^2 - p^{11}) - p^{11}t. \tag{7.4}$$

Substituting these values, and $u = p^{11}$ into the penultimate system of equations returns

$$\begin{aligned} u + 1 + 691v &= t, \\ u^2 + u + 1 + 691x &= t^2 - u, \\ u^3 + u^2 + u + 1 + 691y &= t(t^2 - u) - tu. \end{aligned} \tag{7.5}$$

Eliminating $u$ and $v$, by algebraic means, Grobner basis, or resultant calculation, etc, yields

$$477481x^3 - (2764 + 2073t)t \cdot x^2 + (4 + 4t^2 + 4t^3 + 3t^4) \cdot x + 2764xy - 2(2 + 2t + t^2)t \cdot y + 691y^2 = 0. \tag{7.6}$$

The change of variables $(x, y) \to (-X/691, Y/691)$ in (7.6) produces an elliptic curve in Weierstrass form

$$Y^2 - \frac{4}{691^2}XY - \frac{2t(t^2 - 2t - 2)}{691^6}Y = X^3 + \frac{2073t^2 - 2764t}{691^5}X^2 + \frac{3t^4 + 4t^3 + 4t^2 + 4}{691^8}X \tag{7.7}$$

over the functions field $\mathbb{Q}(t)$. The discriminant and the $j$-invariant are complicated rational functions

$$\Delta(E_t) = r_1(t)/s_1(t), \quad \text{and} \quad j(E_t) = r_2(t)/s_2(t), \tag{7.8}$$

respectively, with $r_i(t), s_i(t) \in \mathbb{Z}[t]$. These are nonvanishing parameters for all but a finite number of $t \in \mathbb{Q}$ such that $r_1(t) = 0$. Therefore, the algebraic curve (7.7) is nonsingular over the rational numbers $\mathbb{Q}$ for any fix integer $t \in \mathbb{Z}$, but a finite number of exceptions.

For each fixed $t \in \mathbb{Q}$, such that $\Delta(E_t) \neq 0$, the elliptic curve (7.7) has a finite number of integral solutions $(x(t), y(t))$, this follows from Theorem 3.8.

Let $H = \max\{|a_i| : 0 \leq i \leq 6\}$ be the height of the elliptic curve (7.7). Now, take a sufficiently large $x \geq H^c$, $c > 2$ constant, by Theorem 3.10, the rank of the curve is bounded by

$$\text{rank}(E(\mathbb{Q})) \leq c_0 \log H^c / \log \log H^c \leq c_1 \log x / \log \log, \tag{7.9}$$

where $c_0, c_1 > 0$ and $c > 2$ are constants.

For each fix $t \in \mathbb{Q}$, let $V_t = \#\{p : \tau_k(p) = t\}$, and $V_t(x) = \#\{p \leq x : \tau_k(p) = t\}$. The quantity $V(t)$ is bounded by the number of integral solutions $\#E_t(\mathbb{Z})$ of the elliptic curve $E_t$ in (7.7). That is, $V(t) \leq \#E_t(\mathbb{Z})$. Therefore, Theorems 3.8 and 3.9 imply that

$$\begin{aligned} V_t(x) \leq V_t &\leq \kappa^{1 + \#S + (1+\delta)\text{rank}(E(K))} \\ &\leq x^{o(1)}. \end{aligned} \tag{7.10}$$

Combining these information yield

$$\#\{\tau_k(n): n \leq x\} \geq \#\{\tau_k(p): p \leq x\}$$
$$\geq \frac{\pi(x) - T}{V_t} \qquad (7.11)$$
$$\geq c_k x^{1+o(1)} / \log x,$$

where $T$ is the number of singular curves, i.e., $\Delta(E_t) = 0$, and $c_k > 0$ is a constant depending on the weight $k = 12, 16, 18, 20, 24, 26$. ∎

It should be observe that $T = 0$ since the system of equations (7.2) to (7.4) can be reparametized, by replacing $p \to 2p$, or $p \to 3p$, etc. in some equations in (7.2) to (7.4) in order to produce a nonsingular algebraic curve.

A list of the values $\tau_k(p)$ for $p < 10^6$, and $k = 12, 16, 18, 20, 22, 26$ was compiled in [GF]. A short version is tabulated below.

| $p$ | $\tau_{12}(p)$ | $\tau_{16}(p)$ | $\tau_{18}(p)$ | $\tau_{20}(p)$ | $\tau_{22}(p)$ | $\tau_{26}(p)$ |
|---|---|---|---|---|---|---|
| 2 | −24 | 216 | −528 | 456 | −288 | −48 |
| 3 | 252 | −3348 | −4284 | 50652 | −128844 | −195804 |
| 5 | 4830 | 52110 | −1025850 | −2377410 | 21640950 | −741989850 |
| 7 | −16744 | 2822456 | 3225992 | −16917544 | −768078808 | 39080597192 |

## 8 The Value $\tau_{12}(n) = 0$

There are several proofs of the congruence $\tau(n) \equiv n\sigma_3(n) \bmod 7$, $\gcd(7, n) = 1$, in [BT], [GL], [RB, Corollary 2.8], and similar references. This congruence serves as a foundation for another proof of the case $k = 12$. Refer to the proof of Theorem 1.2 for other details.

***Theorem 8.1.*** Let $\Delta(z) = \sum_{n \geq 1} \tau(n) q^n \in S_{12}(\Gamma_0(1))$ be a cusp form of weight $k = 12$. Then, the coefficients $\tau(p) \neq 0$ for all primes $p \geq 2$.

***Proof***: Suppose that $\tau(p) = 0$ for some prime $p \geq 11$. Now consider the congruence

$$\tau(p) \equiv p\sigma_3(p) \bmod 7, \qquad (8.1)$$

this is derived from $\tau(n) \equiv n\sigma_3(n) \bmod 7$. Expanding the divisors function $\sigma_s(n) = \sum_{d \mid n} d^s$ produces the Diophantine equations

$$\tau(p) = p(p^3 + 1) + 7v, \qquad (8.2)$$

where $v \in \mathbb{Z}$ are rational integer, and the corresponding congruence equation

$$p(p^3+1)+7v \equiv \tau(p) \bmod 7. \tag{8.3}$$

The hypothesis $\tau(p) = 0$ implies that $\tau(p) \equiv 0 \bmod 7$, and by Theorem 2.2, the subset of primes

$$\{ p \leq x : \tau(p) \equiv 0 \bmod 7 \} \tag{8.4}$$

has positive density. Therefore, equation (8.3) has infinitely many integer solutions. In fact, an explicit sequence of solutions is given by

$$v = -p(p^3+1)/7, \tag{8.5}$$

where the prime $p$ ranges over some arithmetic progressions, for example, $p = 2 \cdot 7m - 1$, $m \geq 1$. This sequence of solutions is the only possibility for both (8.2) and (8.3); other solutions of (8.3) imply that $\tau(p) \neq 0$.

This in turns implies that the congruence $\tau(p) \equiv 0 \bmod 7$ holds over the integers, that is, $\tau(p) = 0$ for each integer solutions of the form (8.5). In particular, the subset of primes

$$\{ p = 2 \cdot 7m - 1 \leq x : \tau(p) = 0 \}, \tag{8.6}$$

has positive density in the set of primes. But, this contradicts the zero density result for the subset of primes $\{ p \leq x : \tau(p) = 0 \}$, see Theorem 2.3. Specifically,

$$\begin{aligned} c_0 \frac{x}{\log x} &\leq \#\{ p = 2 \cdot 7m - 1 \leq x : \tau(p) = 0 \} \\ &\leq \#\{ p \leq x : \tau(p) = 0 \} \\ &\leq c_1 \frac{x}{(\log x)^{3/2-\varepsilon}}, \end{aligned} \tag{8.7}$$

where $c_0, c_1 > 0$ are constants, and $\varepsilon > 0$ is arbitrarily small. Ergo, $\tau(p) \neq 0$ for all primes $p \geq 2$. ∎

As the prime $p \geq 2$ varies over the set of primes, every solution of the equation (8.2) is an integer in the range

$$\frac{-p(p^3+1)-2p^{11/2}}{7} < v < \frac{-p(p^3+1)+2p^{11/2}}{7}, \qquad v \neq \frac{-p(p^3+1)}{7}, \tag{8.7}$$

This follows from the inequality $|\tau(p)| \leq 2p^{(k-1)/2}$, see Theorem 2.5.

**Acknowledgments:** I would like to thank Oliver Rozier for the helpful comments.

## 9. Problems

1. How does the function $p_n \to \tau_k(p_n)$ transforms the prime gaps $p_{n+1} - p_n = O(p_n^{.525})$ to composite gaps $\tau_k(p_{n+1}) - \tau_k(p_n) = O(p_n^{.525(k-1)/2})$?

2. How does the function $p_n \to \tau_k(p_n)$ transforms the prime powers gap $p_{n+1}^2 - p_n^2 = O(p_n^{1.525})$ to prime gaps $\tau_k(p_{n+1}^2) - \tau_k(p_n^2) = O(p_n^{1.525(k-1)/2})$?

3. What is the maximal prime gap $\tau_k(p_{n+m}^2) - \tau_k(p_n^2) \stackrel{?}{=} O(p_n^{1.525(k-1)/2})$?

4. Prove that the asymptotic formula for the primes Fourier coefficients counting function should be of the form

$$\pi_{\tau_k}(x) = \#\{\tau_k(n) \leq x\} = \frac{c_k \, x^{1/(k-1)}}{\log^2 x} + o(\frac{x^{1/(k-1)}}{\log^2 x}), \tag{1}$$

where $c_k > 0$ is a constant. Both the determination of the constant, and the asymptotic order are difficult problems.

5. What is the density of the numbers $\tau_k(n) = a^2 + b^2$, $\tau_k(n) = a^3 + b^3$?

6. Does the function $n \to \tau_k(n)$ has an extention to numbers fields. For example, If $\tau_{12}(2) = -24$, then $\tau_{12}(\alpha) = -2^3$, $\tau_{12}(\beta) = 3$ in some numbers field extention $F$ of the rational numbers $\mathbb{Q}$ such that $2 = \alpha\beta$, with $\alpha \neq 1$, $\beta \neq 1$.

7. For a fixed prime $p$, what are the properties of the orbit $\Theta_p = \{\tau_k(p^{2n}) : n \geq 1\}$? Which orbits $\Theta_p$, $p \geq 1$ contain finitely many primes?

Note: It is known that $\Theta_p = \{\tau_k(p^{2n+1}) : n \geq 1\} = \varnothing$, and $\Theta_p = \{\tau_k(p^n) : n \geq 1\} = \varnothing$, for nonordinary primes.

8. The signs changes of the $p$th coefficients of a cusp forms of large weight $k \geq 12$ for primes in arithmetic progressions. What are the asymptotic formulae:

$$T_+(x) = \#\{\tau(p) > 0 : p \leq x, \, p \equiv a \bmod q\} \stackrel{?}{=} \pi(x)/2\varphi(q), \tag{4.9}$$

and

$$T_-(x) = \#\{\tau(p) < 0 : p \leq x, \, p \equiv a \bmod q\} \stackrel{?}{=} \pi(x)/2\varphi(q),$$

for $0 < a < q$, such that $\gcd(a,q) = 1$ correct?

9. Use the given information: the number $N_p = p + 1 + 2\sqrt{p}\cos\theta_p$ is an independent random variable with respect to the Sato-Tate measure $d\mu = 2\pi^{-1}\sin(\theta)^2 d\theta$, (this is unconditional), and Birch-Swinnerton-Dyer conjecture

$$\prod_{p \leq X} \frac{N_p}{p} = \prod_{p \leq X}\left(1 + \frac{1}{p} + \frac{2\cos\theta_p}{\sqrt{p}}\right) \sim c(\log X)^{\text{rank}(E)}, \tag{2}$$

to complete the evaluation or to estimate the equivalent asymptotic

$$\frac{1}{\log(c\log X)} \sum_{p \leq X} \log\left(1 + \frac{1}{p} + \frac{2\cos\theta_p}{p^{1/2}}\right) = \frac{1}{\log(c\log X)}\left(\sum_{p \leq X}\frac{1}{p} + 2\sum_{p \leq X}\frac{\cos\theta_p}{p^{1/2}} + E(X)\right)$$

$$= 1 + \frac{2}{\log(c\log X)} \sum_{p \leq X} \frac{\cos\theta_p}{p^{1/2}} + O(1) \tag{3}$$

$$\sim \text{rank}(E),$$

where $E(X)$ is some error term is

$$E(X) = O\left(\sum_{p \leq X}\left(\frac{1}{p} + \frac{2\cos\theta_p}{p^{1/2}}\right)^2\right) = O(\log\log x). \tag{4}$$

The difficulty arises in evaluating or estimating

$$R(x) = \sum_{p \leq X} \frac{\cos\theta_p}{p^{1/2}}. \tag{5}$$

If $R(x)/\log\log x \to \infty$ as $x \to \infty$, then the rank is unbounded, otherwise it is bounded.